\documentclass[11pt]{amsart}
\usepackage{amsmath, amssymb, amsbsy, amsfonts, amsthm, latexsym, amsopn, amstext, amsxtra, euscript, amscd, color, mathrsfs}
\usepackage{multirow, hyperref}

\usepackage{todonotes}
\usepackage{lipsum}
\usepackage{marginnote}

\makeatother

\newtheorem{thm}{Theorem}[section]

\newtheorem{thm-con}[thm]{Theorem-Conjecture}
\numberwithin{equation}{section}
\theoremstyle{definition}

\def\cA{{\mathcal A}}
\def\cG{{\mathcal G}} 
\newcommand{\F}{\mathbb F} 
\newcommand{\Tr}{{\operatorname{Tr}}}

\begin{document}
\title[ Classification of Differentially $3$-Uniform Power Mappings]{Toward a Classification of Differentially $3$-Uniform Power Mappings over $\F_{p^2}$ and $\F_{p^3}$}

\author[M. Pal]{Mohit Pal} 
\address{The Selmer Center, Department of Informatics, University of Bergen, P.O. Box 7800, N-5020 Bergen, Norway}
\email{mathmohit@outlook.com}

\begin{abstract}
We study the classification, up to the CCZ-equivalence, of differentially $3$-uniform power mappings over quadratic and cubic extensions of prime fields. Over $\F_{p^2}$ exhaustive searches for primes $3<p<300$ found only the known families E1 and E4 for $p>23$. Over $\F_{p^3}$, the previously unexplained exponents for $7<p<200$ belong to a single family. We give an elementary proof that this family has differential uniformity at most $3$ for primes $p\equiv 3 \pmod 4$, and prove that the bound is attained for $p>1444$.
\end{abstract}
\maketitle
 
\section{Introduction}
Let $\F_{p^n}$ be the finite field with $p^n$ elements, where $p$ is a prime number and $n$ is a positive integer. We denote by $\F_{p^n}^*$, the multiplicative cyclic group of non-zero elements of $\F_{p^n}$. It is well-known that every mapping from $\F_{p^n}$ into itself can be uniquely represented by a polynomial in $\F_{p^n}[x]$ of degree less than $p^n$. Accordingly, we may use the terms function and polynomial interchangeably. A polynomial $f\in \F_{p^n}[x]$ is called a permutation polynomial if the associated mapping $c \mapsto f(c)$ permutes the elements of $\F_{p^n}$. For any function $f: \F_{p^n} \rightarrow \F_{p^n}$, its difference distribution table (DDT) entry at point $(a,b) \in  \F_{p^n} \times \F_{p^n}$, denoted by $\Delta_f(a,b)$, is given by the number of solutions of the equation $f(x+a)-f(x)=b$. The differential uniformity of $f$, denoted by $\Delta_f$, is then given by $\Delta_f= \max \{ \Delta_f(a,b) \mid a \in \F_{p^n}^*, b \in \F_{p^n} \}$. When used as a substitution box in a block cipher, the differential uniformity of a function $f$ quantifies its resistance against the differential attack (see~\cite{KN93}). Lower the differential uniformity, higher is the immunity of the function against differential attacks. 

In the study of the differential uniformity of functions over finite fields, we often classify them with respect to some equivalence relations which preserve this property. Two functions $f, g : \F_{p^n} \rightarrow  \F_{p^n}$ are called linear (affine) equivalent if there exist linear (affine) permutations $A_1, A_2 :  \F_{p^n} \rightarrow  \F_{p^n}$ such that $g = A_2 \circ f \circ A_1$. We say that $f$ and $g$ are extended affine (EA) equivalent if there exist affine permutations $A_1, A_2 :  \F_{p^n} \rightarrow  \F_{p^n}$ and an affine function $A:  \F_{p^n} \rightarrow  \F_{p^n}$ such that $g = A_2 \circ f \circ A_1 +A$. The most general equivalence relation, known so far, which preserves the differential uniformity is the Carlet-Charpin-Zinoviev (CCZ) equivalence~\cite{CCZ1998}.  Two functions $f$ and $g$ are called CCZ-equivalent if there exists an affine permutation $\cA:  \F_{p^n} \times  \F_{p^n} \rightarrow  \F_{p^n} \times  \F_{p^n}$ which maps the graph $\cG_f := (x, f(x))$ to the graph $\cG_g := (x, g(x))$. Thus, we can classify functions over finite fields in CCZ-equivalence classes and then each CCZ-equivalence class can be further classified into EA-equivalent classes. Consequently, the CCZ-class of a function $f$ always contains the EA-class of the function $f$. It is well-known~\cite{BCP06} that if $f$ is a permutation then its CCZ-class also contains the EA-class of $f^{-1}$, the compositional inverse of the function $f$. In the case of power mappings Dempwolff~\cite{D2018} showed that two power mappings $x^{\ell_1}$ and $x^{\ell_2}$ over $\F_{p^n}$ are CCZ-equivalent, if and only if there exists a positive integer $0\leq m < n$, such that $\ell_1 \equiv p^m \ell_2 \pmod{p^n-1}$ or $\ell_1  \ell_2 \equiv p^m \pmod{p^n-1}$.

The lowest possible differential uniformity of a function is $1$ and in this case we say that the function is perfect nonlinear. Perfect nonlinear functions are commonly known as planar functions, and were first introduced by Dembowski and Ostrom~\cite{DO1968} in connection to the study of projective planes. The complete classification of planar functions looks hopeless and is known only over prime fields~\cite{G1990, H1989, J1987, RS1989}. For monomials, the complete classification is known over fields of order $p^2$ due to Coulter~\cite{C2006}, over fields of order $p^3$ due to Bergman et al.~\cite{BCV2022} and over fields of order $p^4$ by Coulter and Lazebnik~\cite{CL2012}. Recently, Beierle and Beyne~\cite{BB2025} completely classified planar monomials for all $n=2^k$ with $p>5$ and $k$ a non-negative integer.

It is well-known that planar functions can never be permutations. Therefore, $2$ is the optimal differential uniformity for permutations over finite fields of odd characteristic and such functions are known as almost perfect nonlinear (APN). The complete classification of APN functions is not known even in the case of monomials. To the best of our knowledge, a systematic study of APN functions in odd characteristic starts with the seminal work of Helleseth, Rong and Sandberg~\cite{HRS1999} and the only known (for $p\geq 7$) infinite classes of APN power maps are those given in~\cite{HRS1999}. We refer interested reader to~\cite{KOPS2026} for detailed analysis on APN power maps in odd characteristic.

Unlike APN functions, very little is known about power mappings over finite fields of odd characteristic with differential uniformity equal to $3$. To the best of our knowledge, there are $5$ known classes (see Table~\ref{Tb1}) of differentially $3$-uniform power maps over finite fields of odd characteristic. In this paper, we present a new infinite class of differentially $3$-uniform power maps over fields of order $p^3$. Experimental results suggest that our family of power maps together with the classes in Table~\ref{Tb1} covers all differentially $3$-uniform power maps over $p^3$ for $7<p<200$.

The paper is organized as follows. In Section~\ref{S2}, we review the state of the art results concerning deifferentially $3$ uniform power maps. In Sections~\ref{S3}, we present a new class of differentially $3$-uniform power maps over $\F_{p^3}$, where $p>1444$.

\section{Preliminaries} \label{S2}
Functions over finite fields with differential uniformity $3$ have been studied for more than three decades, and so far, only five classes of differentially $3$-uniform power maps over finite fields of odd characteristic are known (listed in Table \ref{Tb1}).
{\small 
\begin{table}[h!]
\begin{center}
\begin{tabular}{ |c|c|c|c|c| } 
 \hline
 & $d$ & $p$  & Conditions & Ref\\ 
 \hline
 E1 & $4$ &    &  $p^n>7$ & Folklore\\ 
 \hline
E2 & $p^n-2$ & $3$ &  &  \cite[Remark]{HRS1999}\\
 \hline
 E3 &$\frac{p^n-3}{2}$ & $p \equiv 11,19 \pmod{20}$  &$n$ odd & \cite[Theorem 2]{HS1997}  \\ 
 \hline
 E4 & $\frac{p^n+3}{2}$ & $p\neq 3$  &$p^n \equiv 1 \pmod 4$ & \cite[Theorem 3]{HS1997}  \\ 
 \hline
 E5 & $\frac{p^n-1}{2}+p^{\frac{n+1}{2}}+1$& $p \neq 3$ and $p\neq 17$ when $n=1$ & $p^n \equiv 1 \pmod 4$, $n$ odd & \cite[Theorem 2.1]{Felke2020} \\
 \hline
 \end{tabular}
 \end{center}
 \caption{Known classes of $3$-uniform power maps in odd characteristic.}
  \label{Tb1}
\end{table}
}

It is easy to observe that in the case of prime fields $\F_p$, exponents E4 and E5 coincide. Also, there are no differentially 3-uniform monomials over $\F_3$. Therefore, over prime fields $\F_p$, we have three infinite classes of differentially $3$-uniform power maps, namely, E1, E3 and E5. We performed an exhaustive search for differentially $3$-uniform power maps over prime fields $\F_p$ for primes $p<100$ and the result is listed in  Table~\ref{Tb2}. We could not deduce any new infinite class of differentially $3$-uniform power maps over $\F_p$ from Table~\ref{Tb2}. Due to their vast structural diversity, a complete classification of differentially $3$-uniform power maps over $\F_p$ appears currently to be infeasible.

However, the case in the fields of order $p^2$ seems relatively simpler. Notice that, for $p>3$, we have only two infinite classes, namely, E1 and E4. We performed an exhaustive search for differentially $3$-uniform power maps over $\F_{p^2}$ for primes $3<p<300$. For primes $7 \leq p \leq 23$, we have some unexplained exponents listed in Table~\ref{Tb3}. For every prime $23<p<300$, the differentially $3$-uniform power mappings over $\F_{p^2}$ belongs to E1 or E4. A natural question is whether every differentially $3$-uniform power mappings over $\F_{p^2}$, $p>23$ other than those belonging to E1 and E4. An affirmative answer, together with the exhaustive search classification for $p\leq 23$ would give a complete classification of differentially $3$-uniform power mappings over $\F_{p^2}$.

In the case of fields of order $p^3$, $p>3$, we have four infinite classes of differentially $3$-uniform power maps, namely, E1, E3, E4 and E5. We performed an exhaustive search for differentially $3$-uniform power maps over $\F_{p^3}$ for primes $p<200$ and the result is listed in  Table~\ref{Tb4}. We were trying to get a formula for the unexplained entries in the Table~\ref{Tb4}. Notice that, for $p>7$, the unexplained entries are occurring in the case of $p^n \equiv 3 \pmod 4$. We further divided this into two parts, namely, $p \equiv 3 \pmod 8$ and $p \equiv 7 \pmod 8$, respectively. Our observation is that the unexplained entries belong to the following class:
\begin{equation}\label{expo6}
d=A(p^2+p+1)-1,~\mbox{where}~
A=
\begin{cases}
\frac{p+1}{4}~\mbox{if}~ p \equiv 3 \pmod 8,\\
\frac{3p-1}{4}~\mbox{if}~p\equiv 7 \pmod 8.
\end{cases}
\end{equation} 
In the following section, we shall show that the power mapping $x^d$ over $\F_{p^3}$, where $d$ is an exponent given in Equation~\eqref{expo6}, has differential uniformity less than or equal to $3$.

\section{A new class of power maps over $\F_{p^3}$ with $\Delta_f=3$}  \label{S3}
In this section we present a new class of power mappings with differential uniformity $\leq 3$ over $\F_{p^3}$.

\begin{thm}
Let $p>7$ be a prime with $p \equiv 3 \pmod 4$ and $d=A(p^2+p+1)-1$, where
\begin{equation*}
A=
\begin{cases}
\frac{p+1}{4}~\mbox{if}~ p \equiv 3 \pmod 8,\\
\frac{3p-1}{4}~\mbox{if}~p\equiv 7 \pmod 8.
\end{cases}
\end{equation*} 
Then the differential uniformity of the power map $f(x)=x^d$ over $\F_{p^3}$ is at most $3$. It is equal to $3$ for $p>1444$.
\end{thm}
\begin{proof}
Recall that for any $a,b \in \F_{p^3}$, $a \neq0$,  the difference distribution table (DDT) entry of $f$ at $(a,b) \in \F_{p^3} \times \F_{p^3}$, denote by $\Delta_f(a,b)$, is given by the number of solutions of the following equation
\begin{align}\label{DDTab}
   f(x+a)-f(x) =b \implies  (x+a)^d-x^d=b \implies  \left(\frac{x}{a}+1\right)^d-\left( \frac{x}{a} \right)^d=\frac{b}{a^d}.
\end{align}
Thus, $\Delta_f(a,b)$ is equal to $\Delta_f(1, b/a^d)$. Therefore, in order to compute the differential uniformity of power maps, it is sufficient to compute DDT entries $\Delta_f(1,b)$, $b \in \F_{p^3}$. Now, replacing $x/a$ with $x$, Equation~\eqref{DDTab} reduces to
\begin{align}\label{DDTb}
    (x+1)^d=x^d+b.
\end{align}
Our aim is to show that Equation~\eqref{DDTb} has at most $3$ solutions for all $b \in \F_{p^3}$. It is easy to observe that
\begin{align*}
    (p^2+p)d &=
\begin{cases}
\frac{(p^2+p)(p+1)}{4}(p^2+p+1)-p^2-p~\mbox{if}~ p \equiv 3 \pmod 8,\\
\frac{(p^2+p)(3p-1)}{4}(p^2+p+1)-p^2-p~\mbox{if}~p\equiv 7 \pmod 8,\\
\end{cases}
\equiv \frac{p^3+1}{2} \pmod{p^3-1}.
\end{align*}
Therefore, we have
\begin{align*}
    x^{(p^2+p)d}= \eta(x)x =
    \begin{cases}
        0~&~\mbox{if}~x=0,\\
        x~&~\mbox{if}~x~\mbox{is a nonzero square in}~\F_{p^3},\\
        -x~&~\mbox{if}~x~\mbox{is not a square in}~\F_{p^3},
    \end{cases}
\end{align*}
where $\eta$ is the quadratic character of $\F_{p^3}$. Notice that, since $p \equiv 3 \pmod 4$, we have $p^3 \equiv 3 \pmod 4$ and hence $\eta(-1)=-1$. Also, since $d$ is even and $\gcd(d, p^3-1)=2$, any nonzero square element $u$ in $\F_{p^n}$ has precisely two preimages $u^{p^2+p}$ and $-u^{p^2+p}$ under the map $x^d$. Now, we shall consider three cases, namely, $b=0$, $\eta(b)=-1$ and $\eta(b)=1$.

\textbf{Case 1.} Let $b=0$. In this case, Equation~\eqref{DDTb}, reduces to 
\begin{equation}\label{Eq:c1}
    (x+1)^d=x^d.
\end{equation}
It is easy to see that neither $x=0$ nor $x=-1$ satisfies Equation~\eqref{Eq:c1}. For $x \in \F_{p^3} \backslash \{0,-1\}$, if $x^d=u$ for some nonzero square $u \in \F_{p^n}$, then $x = \pm u^{p^2+p}$. In this case, Equation~\eqref{Eq:c1} reduces to $(x+1)^d=u$ and hence $x+1= \pm u^{p^2+p}$. Thus, the only possible cases are $u^{p^2+p} =1- u^{p^2+p}$ and  $-u^{p^2+p} =1+ u^{p^2+p}$. That is $u^{p^2+p}= \pm \frac{1}{2}$ and since $\eta(-1)=-1$ either $1/2$ or $-1/2$ will be a square in $\F_{p^3}$. It is easy to verify that $x= -\frac{1}{2}$ is a solution of Equation~\eqref{Eq:c1} and hence $\Delta_f(1,0)=1$.

\textbf{Case 2.} Let $\eta(b)=-1$. In this case, RHS of Equation~\eqref{DDTb} is $0$ if and only if $x^d+b=0$. In this case, Equation~\eqref{DDTb} reduces to $(x+1)^d=0 \implies x=-1$ and hence $b=-1$. Now, for any $b \in \F_{p^3} \backslash \{-1\}$ with $\eta(b)=-1$, we shall count the number of solutions $x \in \F_{p^3}$. It is easy to observe that for any $b \in \F_{p^3} \backslash \{-1\}$ with $\eta(b)=-1$, $x \in \{-b^{p^2+p}, b^{p^2+p}\}$ is not a solution of Equation~\eqref{DDTb}. Thus, for any $b \in \F_{p^3} \backslash \{-1\}$ with $\eta(b)=-1$, both LHS and RHS of Equation~\eqref{DDTb} are non-zero.  Assume that $b+x^d=u \implies x= \pm (u-b)^{p^2+p}$ for some $u \in \F_{p^3}^*$ such that $\eta(u-b)=1$. Then, from Equation~\eqref{DDTb}, we have $x+1= \pm u^{p^2+p}$ and $ \eta(u)=1$.  Thus, we have following four subcases:

\textbf{Subcase 2.1.} Let $(u-b)^{p^2+p} = -1+u^{p^2+p}$. Then we have
\begin{align*}
    u^{p^2}b^p+u^pb^{p^2}=1+b^{p^2+p} &\implies v^{p^2}+v^p = \frac{1+b^{p^2+p}}{b^{p^2+p}}\\
     &\implies v^p+v = \frac{1+b^{p+1}}{b^{p+1}}\\
    &\implies v = \frac{Tr(b')-2b'^p}{2},
\end{align*}
where $u=bv$ and $b'= \frac{1+b^{p+1}}{b^{p+1}}$. Thus, we have at most one solution $x$ in this subcase.

\textbf{Subcase 2.2.} Let $(u-b)^{p^2+p} = -1-u^{p^2+p}$. Then we have
\begin{align*}
    -2u^{p^2+p}+u^{p^2}b^p+u^pb^{p^2}=1+b^{p^2+p} 
    &\implies -2v^{p^2+p}+v^{p^2}+v^p= \frac{1+b^{p^2+p}}{b^{p^2+p}}\\
    &\implies  -2v^{p+1}+v^p+v = \frac{1+b^{p+1}}{b^{p+1}}\\
    &\implies  4v^{p+1}-2v^p-2v+1 = \frac{-2-2b^{p+1}}{b^{p+1}}+1\\
    &\implies  (2v-1)^{p+1} = \frac{-2-b^{p+1}}{b^{p+1}}
\end{align*}
where $u=bv$. If RHS of the above equation is zero then we have a unique solution $v= \frac{1}{2}$. If RHS is non-zero then, as $\gcd(p+1, p^3-1)=2$, the above equation has either no solution or has two solutions $v$ and $1-v$. Notice that $\eta(u)=\eta(bv)=-\eta(v)$ and $\eta(u-b)= \eta(b(v-1))=-\eta(v-1)$. Thus, the condition $\eta(u)=\eta(u-b)=1$ implies that $\eta(v)=\eta(v-1)=-1$. Thus, at most one solution from $v$ or $1-v$ is admissible. Thus,  we can have at most one solution $x$ in this subcase. From here we also conclude that the solution $v= \frac{1}{2}$ is not admissible as $\eta(1/2)\neq\eta(-1/2)$.

\textbf{Subcase 2.3.} Let $(u-b)^{p^2+p} = 1-u^{p^2+p}$. Then we have
\begin{align*}
    -2u^{p^2+p}+u^{p^2}b^p+u^pb^{p^2}=-1+b^{p^2+p} 
    &\implies -2v^{p^2+p}+v^{p^2}+v^p = \frac{-1+b^{p^2+p}}{b^{p^2+p}}\\
    &\implies -2v^{p+1}+v^p+v = \frac{-1+b^{p+1}}{b^{p+1}}\\
    &\implies  (2v-1)^{p+1} = \frac{2-b^{p+1}}{b^{p+1}}
\end{align*}
where $u=bv$. Similar to the previous subcase, if RHS of the above equation is zero then we have a unique solution $v= \frac{1}{2}$, which is not admissible as $\eta(1/2)\neq\eta(-1/2)$. If RHS is non-zero then the above equation has either no solution or has two solutions $v$ and $1-v$. The condition $\eta(u)=\eta(u-b)=1$ implies that $\eta(v)=\eta(v-1)=-1$. Thus, at most one solution from $v$ or $1-v$ is admissible. Thus,  we can have at most one solution $x$ in this subcase.

\textbf{Subcase 2.4.} Let $(u-b)^{p^2+p} = 1+u^{p^2+p}$. Then we have
\begin{align*}
    u^{p^2}b^p+u^pb^{p^2}=-1+b^{p^2+p} &\implies v^{p^2}+v^p = \frac{-1+b^{p^2+p}}{b^{p^2+p}} \\
     &\implies v^p+v = \frac{-1+b^{p+1}}{b^{p+1}}\\
    &\implies v = \frac{Tr(b')-2b'^p}{2},
\end{align*}
where $u=bv$ and $b'= \frac{-1+b^{p+1}}{b^{p+1}}$. Thus, we have at most one solution from this subcase.

It is easy to observe thatb for any fixed $b\neq -1$ with $\eta(b)=-1$, if $v$ is a solution in subcase 2.4 then $1-v$ will be a soultion in the subcase 2.1. The condition $\eta(v)=\eta(v-1)=-1$ implies that only one of these can be admissible. So Subcase 2.1 and Subcase 2.1 can contribute at most one solution.  

When $b=-1$, then we have a solution $x=-1$. For $b=-1$, Subcase 2.1 reduces to $v^p+v=2$ which has a unique solution $v=1$ and this does not satisfy the condition $\eta(v)=-1$. Similarly, Subcase 2.4 reduces to $v^p+v=0$ which has a unique solution $v=0$ and it does not satisfy the condition $\eta(v)=-1$. Subcase 2.3 reduces to $(2v-1)^{p+1}=1$, thus $v=0$ or $v=1$ and both fail the character condition. Similarly,  Subcase 2.2 reduces to 
\[
(2v-1)^{p+1}=-3 \implies (2v-1)^p=\frac{-3}{2v-1} \implies (2v-1)^{p^2}=\frac{-3}{2v-1} \implies 2v-1 \in \F_p.
\]
Thus, we have $(2v-1)^{p+1}=-3 \implies (2v-1)^2=-3 \implies v(v-1)= -1$, Which is not possbile as $\eta(v)=\eta(v-1)=-1$. Thus $\Delta(1,-1)=1$.

\textbf{Case 3.} Let $\eta(b)=1$. Since $d$ is even, we have
\[
(-x)^d-(-x-1)^d= x^d-(x+1)^d =-((x+1)^d-x^d).
\]
Thus the involution $x \mapsto -x-1$ gives a bijection between the solutions for $b$ and those for $-b$. Since $\eta(b)=1$ and $\eta(-1)=-1$, we have $\eta(-b)=-1$. Hence, Case 2 applies to $-b$. Therefore, 
\[
\Delta_f(1,b)=\Delta_f(1,-b)\leq 3.
\]

This proves the upper bound.
\end{proof}

\subsection{Attainment of the upper bound}
Let $p>1444$ be a prime with $p\equiv3\pmod4$ and let $d$ be given by Equation~\eqref{expo6}.
The proof of Theorem~3.1 gives $\Delta_f(1,b)\leq3$ for every $b\in\F_{p^3}$.
To prove attainment, it suffices to find a $b$ for which Equation~\eqref{DDTb}
has three distinct solutions. We first give a sufficient condition on $b$
and then show how to obtain such a value.

\par\medskip\noindent\textbf{A sufficient condition.}
Write \(m=p^2+p\), and let \(\eta\) be the quadratic character of
\(\F_{p^3}\). Suppose \(b\in\F_{p^3}\) satisfies
\begin{equation}
\Tr(b)=-1,\qquad \operatorname{N}(b)=1,\qquad
\eta(b+1)=\eta(2b-1)=\eta(-2b-1)=1.
\label{eq:attainment-condition}
\end{equation}
Then \(b^m=b^{-1}\) and \(b\) is a square, since the norm preserves
quadratic character on a cubic extension. By the preimage formula in the
proof of Theorem~3.1, a nonzero square \(u\) has preimages \(\pm u^m\) under \(f\),
where $u^m$ is a square and $-u^m$ is a nonsquare.

For the first solution take \(x=-b^{-2}\).
Indeed, \((b+1)^m=b^{-1}-b\) by the trace and norm conditions, so
\[
x=-(b^2)^m,\qquad x+1=-\bigl(b(b+1)\bigr)^m.
\]
The two quantities inside the powers are nonzero squares and differ
by \(b\). Thus \(f(x+1)-f(x)=b\), with character pair \((-1,-1)\).

For the other two solutions let \(\varepsilon\in\{1,-1\}\).
We seek \(x=-\varepsilon(u-b)^m\) and \(x+1=\varepsilon u^m\), with
\(u,u-b\) nonzero squares. On writing \(h=2u/b-1\), the equation
\(u^m+(u-b)^m=\varepsilon\) becomes
\[
h^m=r_\varepsilon,\qquad r_\varepsilon=2\varepsilon b-1.
\]
This follows from \(u^m+(u-b)^m=(h^m+1)/(2b)\).
By~\eqref{eq:attainment-condition}, \(r_\varepsilon\) is a nonzero square, so the equation has
two opposite roots. Since \(m(m-1)\equiv2\pmod{p^3-1}\), we have
\(h^2=r_\varepsilon^{m-1}\). Also
\(r_\varepsilon^m=4/b+2\varepsilon(b+1)+1\). Hence
\[
u(u-b)=\frac{b^2}{4}(h^2-1)
=\frac{b}{r_\varepsilon}
 \left(1+\frac{1+\varepsilon}{2}b\right).
\]
This is a nonzero square for either sign. The two roots give the
pairs \((u,u-b)\) and \((-(u-b),-u)\). In each pair the characters
agree, and negating reverses both characters because \(\eta(-1)=-1\).
Exactly one root therefore yields two nonzero squares. It gives a
solution with character pair \((-\varepsilon,\varepsilon)\).
The three character pairs are distinct, proving \(\Delta_f(1,b)=3\).

\par\medskip\noindent\textbf{Existence of a suitable derivative value.}
Put \(P_s(T)=T^3+T^2+sT-1\), and set
\[
L_0(s)=s+1,\qquad L_1(s)=5-4s,\qquad L_2(s)=-7-4s.
\]
If \(P_s\) is irreducible over \(\F_p\), any root \(b\) has trace
\(-1\) and norm \(1\). The norms of \(b+1,2b-1,-2b-1\) are,
respectively, \(L_0(s),L_1(s),L_2(s)\), by evaluating \(P_s\) at
\(-1,1/2,-1/2\). Thus a parameter \(s\) for which \(P_s\) is
irreducible and all three \(L_i(s)\) are nonzero squares gives~\eqref{eq:attainment-condition}.

In this subsection write \(\chi\) for the quadratic character of
\(\F_p\), and let \(N_p\) count the parameters just described.
Set
\[
Q(s)=4s^2-5s+23,\quad W(s)=\prod_{i=0}^2(1+\chi(L_i(s))),
\quad \rho(s)=\#\{t\in\F_p:P_s(t)=0\}.
\]
The discriminant is \(D(s)=-L_0(s)Q(s)\). For \(p>23\),
\(2+\chi(D(s))-\rho(s)\) equals \(3\) for an irreducible cubic
and \(0\) otherwise. For separable cubics this follows from the
sign of Frobenius on the roots. A repeated-root cubic has two
distinct rational roots, since a triple root would force \(28=0\).
At a zero of any \(L_i\), \(P_s\) has a rational root. Consequently
\[
24N_p=\sum_{s\in\F_p}W(s)(2+\chi(D(s))-\rho(s)).
\]

Put
\[
\begin{aligned}
F_0(t)&=t(1-t),\\
F_1(t)&=t(2t-1)(2t^2+3t+4),\\
F_2(t)&=t(2t+1)(2t^2+t-4)
\end{aligned}
\]
\[
A=\sum_s W(s),\quad B=\sum_s\chi(-Q(s))W(s),\quad
C=\sum_t\prod_{i=0}^2(1+\chi(F_i(t))).
\]
Then \(24N_p=2A+B-C+1\). To see this, the discriminant term
equals \(B+\chi(2)W(-1)\), since
\(\chi(L_0)(1+\chi(L_0))=1+\chi(L_0)\) except at \(s=-1\),
where \(\chi(-Q(-1))=-\chi(2)\).
For the root term write \(s=\phi(t)=t^{-1}-t^2-t\), with \(t\ne0\).
The pullbacks \(L_i(\phi(t))\) are \(F_i(t)/t^2\) for \(i=1,2\),
and \((t+1)^2F_0(t)/t^2\) for \(i=0\).
The contribution to \(C\) at \(t=0\) is \(1\). At \(t=-1\) it is
\((1-\chi(2))W(-1)\), whereas the root term there is \(W(-1)\).
Thus the root term equals \(C-1+\chi(2)W(-1)\).

For completeness, the following table records the expansions used
with the Weil bound for quadratic character sums \footnote{See N. M. Katz, \emph{Sato--Tate theorems for finite-field Mellin transforms}, Introduction, p.~3, \url{https://web.math.princeton.edu/~nmk/mellin246.pdf}.}.
\[
\begin{array}{c|l|c|c}
\text{sum}&\text{degrees of nonconstant character arguments}
   &\text{central value}&\text{error bound}\\ \hline
 A&1,1,1,2,2,2,3 &p+1 &2\sqrt p\\
 B&2,3,3,3,4,4,4,5 &0 &16\sqrt p\\
 C&2,4,4,4,4,6,8 &p+2 &18\sqrt p
\end{array}
\]
For a squarefree polynomial of degree \(j\), the sum is bounded by
\((j-1)\sqrt p\) in odd degree. In even degree it differs from
\(-\chi(\text{leading coefficient})\) by at most \((j-2)\sqrt p\).
The arguments for \(A\) and \(B\) use the pairwise coprime
squarefree factors \(L_0,L_1,L_2,Q\). Indeed, \(Q\) has
discriminant \(-343\) and takes values \(32,23,44\) at the
three roots of \(L_0,L_1,L_2\).
For \(C\), cancel the square factor \(t^2\) from the products
of two or three \(F_i\). The resulting arguments are squarefree:
the two quadratic factors have discriminants \(-23,33\), are
nonzero at \(0,1,\pm1/2\), and can share only the root \(-4\)
where both have value \(24\). The leading-coefficient characters
for \(A\)'s quadratic terms are \((-1,-1,1)\).
In \(B\) they are \(-1\) for degree two and \((1,1,-1)\)
for degree four. In \(C\), in the order
\(F_0,F_1,F_2,F_0F_1,F_0F_2,F_1F_2,F_0F_1F_2\)
after cancellation, they are \((-1,1,1,-1,-1,1,-1)\).
The reduced pair terms have characters \((-1,-1,1)\) at \(0\),
so excluding their spurious contribution adds \(1\) to the
central value for \(C\).

It follows that \(24N_p\ge p+1-38\sqrt p>0\) for \(p>1444\).
Thus an admissible parameter exists for every prime $p>1444$ with $p\equiv3\pmod4$. The construction above gives a value $b$ with $\Delta_f(1,b)=3$, proving attainment in this range.

\section*{Declaration on the use of AI tools}

The author used ChatGPT to assist in developing and drafting the proof that the upper bound on differential uniformity is attained. The author reviewed the argument, checked the mathematical reasoning, and takes full responsibility for the final proof and the content of the manuscript.

\newpage
\section*{Appendix. Computational tables}
\begin{table}[h!]
\begin{center}
\begin{tabular}{ |c|c|c|c|l|} 
 \hline
 $p$ & E1 & E3 & E4 & Unexplained  \\ 
 \hline
 $5$ & $4$ &  &  &\\
 \hline
$7$ & &  & &\\
 \hline
 $11$ & $4$ &&& 6 \\
 \hline
 $13$ & $4$&&$8$& \\
 \hline
 $17$ & $4$&&& $ 6, 12, 14$ \\
 \hline
 $19$ & $4$ &$8$&& $16$ \\
 \hline
 $23$ & $4$ &&& $16$\\
 \hline
 $29$ & $4$ &&$16$& $6,8, 12, 18, 20, 22, 24$ \\
 \hline
 $31$ & $4$ &$14$&& $26, 28$ \\
 \hline
 $37$ & $4$ &&$20$& $8, 10, 16, 22, 28, 32$ \\
 \hline
 $41$ & $4$  &&$22$&$6, 12, 14, 18, 26, 28, 34, 36$\\
 \hline
 $43$ & $4$ &&& $ 10, 16, 26, 38, 40$\\
 \hline
 $47$ & $4$ &&& $18, 20, 26, 30, 32, 34, 38$ \\
 \hline
 $53$ &  $4$ & &$28$& $6,8,10,12,16,20,22,30,32,34,36,42,44,46,48$ \\
\hline
 $59$ & $4$ & $28$&& $10,12,14,22,24,32,46,50$ \\
\hline
 $61$ &  $4$ & &$32$& $14,22,26,28,38,44,52,56$ \\
\hline
 $67$ &  $4$ & && $8,14,16,20,26,28,38,40,46,62$\\
\hline
 $71$ &  $4$ & $34$&& $8,12,16,22,26,38,48,52,64,66$\\
\hline
 $73$ &  $4$ & &$38$& $10,14,22,26,28,34,44,46,50,52,58,62,70$ \\
\hline
 $79$ &  $4$ & $38$&& $14,16,32,34,44,46,58, 70$ \\
\hline
 $83$ &  $4$ & && $8,12,16,24,26,34,38,44,46,68,70,72$ \\
\hline
$89$ & $4$ & &$46$& $6,14,20,26,28,38,42,50,52,58,60,62,70,74,76,78,82$ \\
\hline
 $97$ &  $4$ & &$50$& $10,26,28,46,58,68,70,74,76,82,86$ \\
\hline
 \end{tabular}
 \end{center}
 \caption{Differentially 3-uniform power maps $x^d$ over $\F_p,~ p<100$.}
 \label{Tb2}
\end{table}

\begin{table}[h!]
\begin{center}
\begin{tabular}{ |c|c|c|l|} 
 \hline
 $p$ & E1  & E4 & Unexplained  \\ 
 \hline
 $5$ & $(4, 20)$ & $(14, 22)$ & \\
 \hline
 $7$ & $(4, 28)$ & $(26, 38)$ & $(10, 22)$\\
 \hline
 $11$ & $(4,44)$ &$(62, 82)$ & $(52, 92)$ \\
 \hline
 $13$ & $(4,52)$ &$(86, 110)$ & $ (44, 68), (46, 94), (74, 122)$ \\
 \hline
 $17$ & $(4,68)$ &$(146, 178)$ & $(22, 86), (116, 244), (166, 230), (218, 250)$ \\
 \hline
 $19$ & $(4,76)$ &$(182, 218)$ & $(106, 214)$ \\
 \hline
 $23$ & $(4,92)$ &$(266, 310)$ & $(74, 118), (76, 164), (178, 398)$ \\
 \hline
\end{tabular}
\end{center}
 \caption{Differentially 3-uniform power maps $x^d$ over $\F_{p^2}, 3< p \leq 23$.}
 \label{Tb3}
\end{table}

{\small
\begin{table}[h!]
\centering
\begin{tabular}{ |c|c|c|c|c|l|} 
 \hline
 $p$ & E1 & E3 & E4 & E5 & Unexplained  \\ 
 \hline
 11 & (4,44,484) &(544, 664, 654) & &&(278, 398, 388) \\
 \hline
 13 & (4, 52, 676) & & (1100, 1124, 1436) &  (1112, 1280, 1268) & \\
 \hline
 17 & (4, 68, 1156) & & (2458, 2490, 3034) & (2474, 2762, 2746) & \\
 \hline
 19 & (4, 76, 1444) &(3068, 3428, 3410) & & & (1544, 1904, 1886) \\
 \hline
 23 & (4, 92, 2116)&&&&(8872, 9400, 9378) \\
 \hline
 29 &(4, 116, 3364) & &(12252,13876, 12196) &(12224,13064, 13036) & \\
 \hline
 31 & (4, 124, 3844) & (13934, 14894, 14864) & & & (21878,22838,22808) \\
 \hline
 37 & (4, 148, 5476) & & (25328,25400,28064) & (25364,26732,26696) &  \\
 \hline
\end{tabular}
 \caption{Differentially 3-uniform power maps over $\F_{p^3}$.}
 \label{Tb4}
\end{table}
}
\end{document}